\documentclass[journal]{IEEEtran}
\usepackage{comment}
\usepackage{cite}
\usepackage{hyperref}

\usepackage{multirow}

\ifCLASSINFOpdf
\else
\fi
\usepackage{amsmath}
\usepackage{amsfonts}
\usepackage{algorithmic}
\usepackage[ruled,vlined]{algorithm2e}
\allowdisplaybreaks
\usepackage{color}
\begin{document}
%
\title{An MISOCP-Based Decomposition Approach for the Unit Commitment Problem \\ with AC Power Flows}
%
%
%

\author{Deniz Tuncer, Burak Kocuk
\thanks{
D. Tuncer and B. Kocuk are with the Faculty of Engineering and Natural Sciences, Sabanc{\i} University, Istanbul, Turkey (e-mail: dtuncer, burak.kocuk@sabanciuniv.edu).
}
}

\maketitle

\begin{abstract}
Unit Commitment (UC) and Optimal Power Flow (OPF) are two fundamental problems in short-term electric power systems planning that are traditionally solved sequentially. The state-of-the-art mostly uses a direct current flow approximation of the power flow equations in the UC-level and the generator commitments obtained are sent as input to the OPF-level. However, such an approach can yield infeasible or suboptimal generator schedules. In this paper, we aim to simultaneously solve the UC Problem with alternating current (AC) power flow equations, which combines the challenging nature of both UC and OPF problems. Due to the highly nonconvex nature of the AC flow equations, we utilize the mixed-integer second-order cone programming (MISOCP) relaxation of the UC Problem as the basis of our solution approach. For smaller instances, we develop two different algorithms that exploit the recent advances in the OPF literature and obtain high-quality feasible solutions with provably small optimality gaps. For larger instances, we propose a novel Lagrangian decomposition based approach that yields promising results.
\end{abstract}


\section*{Nomenclature}
\addcontentsline{toc}{section}{Nomenclature}
Decision Variables
\begin{IEEEdescription}[\IEEEusemathlabelsep\IEEEsetlabelwidth{$\overleftarrow{p}_{ij,t}$ ($\overleftarrow{q}_{ij,t}$)}]
\item[$|V_{i,t}|$]                                                  Voltage magnitude of bus $i$ in period $t$
\item[$\theta_{i,t}$]                                               Voltage phase angle of bus $i$ in period $t$
\item[$p_{i,t}^{g}$ ($q_{i,t}^{g} $)]                               Active (reactive) power generation amount of generator $i$ in period $t$
\item[$\overrightarrow{p}_{ij,t}$ ($ \overrightarrow{q}_{ij,t}$)]   Active (reactive) power through line $(i,j)$ in period $t$ in forward direction
\item[$\overleftarrow{p}_{ij,t}$ ($\overleftarrow{q}_{ij,t}$)]      Active (reactive) power through line $(i,j)$ in period $t$ in backward direction
\item[$u_{i,t}$]                                                    Commitment status of generator $i$ in period $t$
\item[$v_{i,t}$]                                                    Startup status of generator $i$ in period $t$
\item[$w_{i,t}$]                                                    Shutdown status of generator $i$ in period $t$
\item[$f_{ij,t}$] Real power flow across the line $(i,j)$.

\item

\end{IEEEdescription}

Parameters
\begin{IEEEdescription}[\IEEEusemathlabelsep\IEEEsetlabelwidth{$f_{i,t}(u_{i,t},v_{i,t},w_{i,t})$}]
\item[$p_{i,t}^{d}$($q_{i,t}^{d}$)] Active (reactive) energy demand of bus $i$ in period $t$
\item[$\delta(i)$ ]                 Set of neighbors of bus $i$
\item[$g_{ii}$($b_{ii}$)]           Shunt conductance (susceptance) of bus $i$ 
\item[$p_i^{\text{min}}$($p_i^{\text{max}}$)]   Minimum (maximum) active generation limit of generator $i$
\item[$q_i^{\text{min}}$($q_i^{\text{max}}$) ]  Minimum (maximum) reactive generation limit of generator $i$
\item[$Y^{ij}$ ]                    Admittance matrix of line $(i,j)$
\item[$ G_{ij}^{ff} $ and $ B_{ij}^{ff}$ ]  Real and imaginary parts of $Y^{ij}_{11}$
\item[$ G_{ij}^{ft} $ and $ B_{ij}^{ft}$ ]  Real and imaginary parts of $Y^{ij}_{12}$
\item[$ G_{ij}^{tf} $ and $ B_{ij}^{tf}$ ]  Real and imaginary parts of $Y^{ij}_{21}$
\item[$ G_{ij}^{tt} $ and $ B_{ij}^{tt}$ ]  Real and imaginary parts of $Y^{ij}_{22}$
\item[$\overline{S}_{ij}$  ]                Maximum power that can pass through line $(i,j)$
\item[$C_{i,t}(p_{i,t}^{g})$  ]             Operational cost of generator $i$ in period $t$
\item[$RU_{i}$ ($RD_{i}$)]                  Ramp up (down) rate of generator $i$
\item[$MinUp_{i}$ ($MinDw_{i}$) ]           Minimum uptime (downtime) of generator $i$
\item[$f_{i,t}(u_{i,t},v_{i,t},w_{i,t})$]   Fixed cost of generator $i$ in period $t$
\item[$\lambda_{i,t}^{up}$] Lagrangian multiplier for minimum uptime constraint \eqref{eq:minupu}
\item[$\lambda_{i,t}^{dw}$] Lagrangian multiplier for minimum downtime constraint \eqref{eq:mindwu}
\item[$\lambda_{i,t}^{log}$] Lagrangian multiplier for the logical constraint \eqref{eq:logicu}
\item[$\lambda_{i,t}^{rd}$] Lagrangian multiplier for the ramp down constraint \eqref{eq:rampdwu}
\item[$\lambda_{i,t}^{ru}$] Lagrangian multiplier for the ramp up constraint \eqref{eq:rampupu}.
\end{IEEEdescription}

%
\IEEEpeerreviewmaketitle

\section{Introduction}
%
%
%
%
\IEEEPARstart{T}{he} Optimal Power Flow (OPF) and the Unit Commitment (UC) Problems are the backbones of the short-term power system planning. The objective of the UC Problem is to determine the most cost-efficient schedule of generators for the next day whereas the objective of the OPF Problem is to find a minimum cost feasible electricity dispatch in a power network given the generator schedules. The UC Problem is typically solved with the approximate direct current (DC) power flow equations whereas the OPF Problem requires the AC power flow equations to be taken into account.
Although  it is common to solve these two problems sequentially, simultaneously  deciding the generator schedules and power dispatch can yield lower total cost. This is the main motivation and the goal of our paper.


Since OPF is a fundamental problem, the related literature is quite rich. There are various methods for solving the OPF Problem, which we categorize into three groups: i)  nonlinear programming (NLP) methods that aim to locate locally optimal solutions, ii) linear programming (LP) methods that approximately solve the problem, iii) conic relaxations based on semidefinite programming (SDP) relaxations and second order cone programming (SOCP) that can provide dual bounds.

The first group of methods used in the OPF literature is the NLP based methods. They generally focus on obtaining the stationary points of the problem, utilizing Newton's method or interior point methods. Newton's method requires using Lagrangian multipliers as penalty terms for the constraints in the constrained optimization problems \cite{frank1}. The papers \cite{sasson} and \cite{dacosta} utilize the Newton's method for solving the OPF problem. By transforming the optimality conditions, \cite{tognola} makes the problem solvable by an algorithm solely based on the Newton-Raphson method. While these methods are good at finding stationary points, they may get stuck at a local optimum point since the OPF problem is nonconvex \cite{waqquas}. Also, since the initial point is important for these methods they may fail to converge to a solution \cite{kocuk016}.

The second group of algorithms in the OPF literature uses  LP based methods, which utilize the DC approximation of the AC power flows \cite{lpapproxcoffrin,stottdc,bienstock}. Since  LPs can be solved efficiently, LP based methods are widely used in industry to solve the OPF Problem. LP methods are also used as a subroutine in a successive linear programming (SLP) framework \cite{griffith}. However, the LP methods may not produce  an AC feasible dispatch since the reactive powers and losses are completely ignored or approximated.

The third group of algorithms in the OPF literature utilizes conic programming relaxations.
In particular, SDP relaxations \cite{baiwei1,baiwei2} have drawn significant interest due to their strength and polynomial solvability. Such relaxations can provide globally optimal solutions when the relaxation is exact, although this can be only guaranteed for very special cases  \cite{lavaei}.  Penalized SDP problems are also proposed to find approximately feasible solutions \cite{madani}.
The computational burden of solving SDP relaxations can be relieved by cheaper but weaker relaxations in the form of the SOCP relaxation \cite{jabrcone}. The authors of \cite{strongsocpkocuk} exploit the computational advantages of the SOCP relaxation and strengthen it using three types of inequalities so that their approach is both more efficient and accurate than the standard SDP relaxation. %
However, despite all the effort, conic programming relaxations are not guaranteed to provide feasible solutions and they are very useful  for dual bounding purposes.

We have mentioned the rich literature on the OPF Problem. However, the literature on the Multiperiod OPF (MOPF) Problem, which can be seen as a generalization of the OPF Problem and a special case of the UC Problem, problem is somewhat limited. A generalized Benders decomposition is used in \cite{alguacil}, in which network constraints are modeled according to a DC approximation. An LP based interior point algorithm is used in \cite{demirovic}, where Newton-Raphson method is used as a corrector. The problem is formulated as a two stage nonlinear program in \cite{schanen}, followed by solving the problem with an interior point method.

The literature on the UC Problem is also quite rich. 
Over the years, many different approaches for the problem are developed which we categorized into three groups: i) Lagrangian relaxation, ii) mixed-integer linear programming (MILP), iii) Benders decomposition, and iv) conic relaxations.

The first method we would like to mention for the UC Problem is the Lagrangian Relaxation (LR) method, which aims to decompose it into smaller subproblems. A branch-and-bound algorithm is proposed in \cite{muckstadt} where they use a Lagrangian method and they decompose the problem into problems with a single generator. LR is utilized in \cite{dubost} in order to find a lower bound and a primal relaxed solution, where the solution is used in a heuristic resolution method. The main drawback of these methods is that the solution obtained by solving LR dual problems might not be primal feasible for the UC Problem \cite{beltran02}.

The second method, mixed integer linear program (MILP) based methods depend on the DC approximation of the power flow equations and is the state-of-the-art for the UC Problem. 
For tightening the MILP formulation, different methods and formulations are proposed in \cite{carroyo,morales13,frangioni09,ostrowski12}.
The drawback of this method is that since it disregards the reactive powers, the solutions might cause stress in the power network. One exception in this direction is \cite{nanou21}, which incorporates AC power flow constraints are incorporated into the MILP model as piecewise linear constraints. 

The third method we would like to mention is the Benders Decomposition (BD) \cite{fuetal,sifuentes}. These references have a master problem and nonlinear subproblems, from which they obtain Benders cuts for the master problem. The drawback of the BD approach is that it converges slowly and is computationally expensive \cite{bendersslow}.

There is only very limited literature on  the UC Problem with AC power flow equations and this literature typically utilize  conic relaxation methods such as SDP relaxation  \cite{baiwei2}. 
A strengthened SDP problem is solved in \cite{atamturk}, which is a convex problem and its minimum can be found efficiently. The SDP relaxation together with Benders' decomposition is considered in \cite{paredes}.
A mixed-integer second-order cone programming (MISOCP) master problem with MOPF subproblems  are considered in \cite{anyacastillo,liuetal}. However, their approaches fail to produce feasible generation schedules within a practical time limit. 
Our paper contributes to this literature by proposing MISOCP based methods that can provide  provably high quality feasible solutions within a reasonable computational budget.

In our study, we consider the UC problem with AC power flows. We aim to find AC feasible commitment schedules for the generators, which requires us to include AC OPF constraints in the UC formulation. Therefore, the problem we are interested in is quite challenging. We utilize the mixed-integer conic relaxation of the UC Problem and solve an MISOCP in order to obtain a candidate generator commitment schedule and also a lower bound for the problem. Then,  given the commitment schedule, we solve an MOPF problem and find an AC feasible solution using a local solver. The lower bounds are compared with the upper bounds obtained from the feasible solutions, and we are able to provide an optimality gap.

For solving the UC problem with AC power flows, we propose three different MISOCP-based algorithms. The first algorithm we propose solves an MISOCP, followed by an NLP problem. In the second algorithm, we strengthen the MISOCP problem by utilizing two methods from \cite{strongsocpkocuk}. These algorithms enable us to solve the UC Problem with AC power flows efficiently, and within 1.3\% optimality gap for small-size instances. However, for large-scale instances solving a 24-period MISOCP problem is quite challenging for which we also develop a novel  decomposition algorithm. This algorithm divides the planning horizon into equally sized blocks, and solves Lagrangian MISOCP subproblems. With the decomposition algorithm, we are able to solve larger instances within a time limit that we set to be one hour, with small optimality gaps. In addition, since there is a lack of publicly available AC UC problem instances, we create and publish our own instances based on publicly available AC OPF instances. 

The rest of the paper is organized as follows: In Section \ref{problemdef}, we define the problems and formulations. In Section \ref{solutionalgos}, we introduce  our proposed solution algorithms. In Section \ref{compresults}, we report the computational results. We conclude our paper in Section \ref{conclusion} with final remarks and comments on possible future work.

\section{Problem Definitions and Formulations}\label{problemdef}
\subsection{Multiperiod Optimal Power Flow}

Consider a power network $\mathcal{N} = (\mathcal{B},\mathcal{L})$ in which the set $\mathcal{B}$ denotes the set of buses and $\mathcal{L}$ denotes the set of transmission lines. Generators are attached to the buses, which are denoted by the set $\mathcal{G} \subseteq \mathcal{B}$. Given the demand of each bus at each time period in a planning horizon $\mathcal{T}$, the MOPF Problem aims to determine the amount of energy produced at each generator, and the economic dispatch of the active and reactive energy. 

Note that the single-period OPF Problem is just a special case of MOPF Problem with $|\mathcal{T}|=1$ and without the ramping constraint.
On the other hand, the MOPF Problem is a special case of the UC Problem for the given generator commitment statuses.

\subsubsection{Formulation}
\label{sec: mopf formulation}
The classical formulation of the OPF and MOPF Problems either uses the rectangular (or $e-f$) formulation  or the polar (or $V- \theta$) formulation. We instead use the alternative formulation first proposed by \cite{exposito} for the OPF Problem, which uses   $c-s-\theta$ variables as follows:
\begin{subequations}\label{eq:csvardeclarations}
\begin{align}
 &\hspace{0.25em} c_{ii,t} := |V_{i,t}|^2 &i \in \mathcal{B},t \in \mathcal{T} \\
  & \hspace{0.25em} c_{ij,t} :=\ |V_{i,t}||V_{j,t}|\cos(\theta_{i,t} - \theta_{j,t}) &(i,j) \in \mathcal{L},t \in \mathcal{T}\\
  & \hspace{0.25em} s_{ij,t} := -|V_{i,t}||V_{j,t}| \sin(\theta_{i,t} - \theta_{j,t}) &(i,j) \in \mathcal{L}, t \in \mathcal{T}
\end{align}
\end{subequations}

We are now ready to provide the formulation of the MOPF Problem:
\begin{subequations}\label{eq:themopf}
\begin{align}
\mathrm{\min} &\hspace{0.25em}  \sum_{t \in \mathcal{T}} \sum_{i \in \mathcal{G}} C_i(p_{i,t}^g)  \label{eq:objective} \\
  \mathrm{s.t.} & \hspace{0.25em} \text{For each }  i \in \mathcal{B}, t \in \mathcal{T}: \notag \\
  &\hspace{0.25em} p_{i,t}^{g} - p_{i,t}^{d} =  \ g_{ii} c_{ii,t} +  \sum_{j\in\delta(i)}( \overrightarrow{p}_{ij,t} + \overleftarrow{p}_{ij,t}) \label{eq:cstheta active balanceu} \\
  & \hspace{0.25em} q_{i,t}^{g} - q_{i,t}^{d} =   -b_{ii} c_{ii,t} +  \sum_{j\in\delta(i)}( \overrightarrow{q}_{ij,t} + \overleftarrow{q}_{ij,t})  \label{eq:cstheta reactive balanceu} \\
  & \hspace{0.25em} q_i^{\text{min}}  \le q_{i,t}^{g} \le q_i^{\text{max}} \label{eq:vtheta reactive gen}\\
  & \hspace{0.25em} p_i^{\text{min}}  \le p_{i,t}^{g} \le p_i^{\text{max}} \label{eq:vtheta active gen} \\
  & \hspace{0.25em} -RD_i \leq p_{i,t}^{g}-p_{i,t'}^{g} \leq RU_i \label{eq:rampopf}\\
  & \hspace{0.25em} \text{For each }  (i,j) \in \mathcal{L}, t \in \mathcal{T}: \notag \\
    & \hspace{0.25em} \overrightarrow{p}_{ij,t} = G_{ij}^{ff}c_{ii,t} + G_{ij}^{ft} c_{ij,t} -B_{ij}^{ft} s_{ij,t}  \label{eq:cstheta real flow fwu} \\
  & \hspace{0.25em} \overleftarrow{p}_{ij,t} = G_{ij}^{tt}c_{jj,t} + G_{ij}^{tf} c_{ij,t} + B_{ij}^{tf} s_{ij,t} \label{eq:cstheta real flow bwu} \\
    & \hspace{0.25em} \overrightarrow{q}_{ij,t} = -B_{ij}^{ff}c_{ii,t} - B_{ij}^{ft} c_{ij,t} -G_{ij}^{ft} s_{ij,t} \label{eq:cstheta reactive flow fwu}\\
  & \hspace{0.25em} \overleftarrow{q}_{ij,t} = -B_{ij}^{tt}c_{jj,t} - B_{ij}^{tf} c_{ij,t} +G_{ij}^{tf} s_{ij,t} \label{eq:cstheta reactive flow bwu} \\
  &\hspace{0.25em} (\overrightarrow{p}_{ij,t})^2 + (\overrightarrow{q}_{ij,t})^2 \le \overline{S}_{ij}^{2}  \label{eq:ubline fw active}\\
  & \hspace{0.25em} (\overleftarrow{p}_{ij,t})^2 + (\overleftarrow{q}_{ij,t})^2 \le \overline{S}_{ij}^{2} \label{eq:ubline bw active}\\
   &\hspace{0.25em} (c_{ij,t})^{2}+(s_{ij,t})^{2} = c_{ii,t} c_{jj,t}  \label{eq:vtheta consistency1} \\
 &\hspace{0.25em} s_{ij,t} = \tan (\theta_{j,t} -\theta_{i,t})  c_{ij,t}.\label{eq:vtheta consistency2}
\end{align}
\end{subequations}

The objective function \eqref{eq:objective} is a quadratic function of the real power generation amount. Constraints \eqref{eq:cstheta active balanceu} and \eqref{eq:cstheta reactive balanceu} are real and reactive power balances at bus $i$, respectively. Constraints \eqref{eq:vtheta reactive gen} and \eqref{eq:vtheta active gen} are bounds for reactive and real power outputs of generator $i$, respectively.  Constraint \eqref{eq:rampopf} is the ramping constraint, which ensures that the active power generation amount does not change rapidly from one period to another. 
Notice that the ramping constraint requires the historical data regarding the real power dispatch variables. To model this aspect of the problem correctly,  we  use a parameter $t'$ in this constraint. If we have the historical data for $t=0$, then $t'$ is simply set to $t-1$. Since  this information is not present in our instances, we assume that  the demand pattern is cyclic and define
\[
t' = 
\begin{cases}
t-1 & \text{if } t \neq 1 \\
24 & \text{if } t = 1
\end{cases}.
\]
Constraints \eqref{eq:cstheta real flow fwu} and \eqref{eq:cstheta real flow bwu} are real power flows from bus $i$ to $j$ and from bus $j$ to $i$, respectively. Constraints \eqref{eq:cstheta reactive flow fwu} and \eqref{eq:cstheta reactive flow bwu}
are reactive power flows from bus $i$ to $j$ and from bus $j$ to $i$, respectively.  Constraints \eqref{eq:ubline fw active} and \eqref{eq:ubline bw active} are power flow upper bounds for line $(i,j)$, for forward and backward flow, respectively.   Constraints \eqref{eq:vtheta consistency1} and \eqref{eq:vtheta consistency2} are the consistency constraints for the relation between $c-s$ variables and original variables of the problem.

\subsubsection{SOCP Relaxation}
\label{sec:socp relax}

An SOCP relaxation of \eqref{eq:themopf} is obtained by relaxing the constraint \eqref{eq:vtheta consistency1} and omitting the constraint \eqref{eq:vtheta consistency2}.
\begin{subequations}\label{eq:socprelx}
\begin{align}
\mathrm{\min} &\hspace{0.25em}  \sum_{t \in \mathcal{T}} \sum_{i \in \mathcal{G}} C_i(p_{i,t}^g)  \label{eq:altobj} \\
  \mathrm{s.t.}   &\hspace{0.25em} (c_{ij,t})^{2}+(s_{ij,t})^{2} \leq c_{ii,t} c_{jj,t} &(i,j) \in \mathcal{L}, t \in \mathcal{T} \label{eq:conicrelax} \\
	&\hspace{0.25em} \eqref{eq:cstheta active balanceu}-\eqref{eq:ubline bw active}. \notag
\end{align}
\end{subequations}

This SOCP relaxation is utilized for constructing an MISOCP relaxation of the UC problem with AC power flows in Section \ref{sec:misocp relax}.

\subsubsection{DC Approach}

Since LP's are efficiently solvable, linearized AC power flow equations are used with the assumption of a DC power grid. The DC approximation of the AC MOPF is formulated as follows:

\begin{subequations}\label{eq:dcapprox}
\begin{align}
\mathrm{\min} &\hspace{0.25em}  \sum_{t \in \mathcal{T}} \sum_{i \in \mathcal{G}} C_i(p_{i,t}^{g})  \label{eq:obj dcopf} \\
  \mathrm{s.t.} &\hspace{0.25em} \text{For each }i \in \mathcal{B}, t \in \mathcal{T}: \notag \\
  &\hspace{0.25em} p_{i,t}^{g} - p_{i,t}^{d} = \sum_{j \in \delta(i)} f_{ij,t}  \label{eq:dcbalance} \\
  &\hspace{0.25em} \text{For each }(i,j) \in \mathcal{L}, t \in \mathcal{T}: \notag \\
  & \hspace{0.25em} f_{ij,t} = -f_{ji,t} = B_{ij}(\theta_{i,t}-\theta_{j,t} ) \label{eq:dcthetaconst}\\
  & \hspace{0.25em} -\overline{S}_{ij} \leq f_{ij,t} \leq \overline{S}_{ij} \label{eq:dcsecurity}\\
  & \hspace{0.25em} \eqref{eq:vtheta active gen}, \eqref{eq:rampopf}. \notag
\end{align}
\end{subequations}

In this formulation, constraint \eqref{eq:dcbalance} enforces flow conservation in each bus $i$. Constraint \eqref{eq:dcthetaconst} defines real power across the line $(i,j)$. Constraint \eqref{eq:dcsecurity} enforces an upper bound on the power over line $(i,j)$. We will make use of this formulation to construct an MILP approximation of the UC problem with AC power flows in Section \ref{sec:milp approx}.

\subsection{Unit Commitment}
The aim of the UC Problem is to determine the commitment schedule of generators over a time horizon given the energy demand of the power system in a way that  minimizes system-wide energy generation cost, fixed and startup costs of generators  while satisfying the operational requirements. Generally, the commitment schedule is obtained with respect to an approximate DC power flow in the UC level. In our study, we aim to solve the UC  Problem while considering the AC power flows, therefore, the problem is challenging.

\subsubsection{Formulation}

We present our UC formulation with AC power flows below that utilizes the MOPF formulation from Section \ref{sec: mopf formulation} with the addition of commitment, startup and shutdown status variables. Recall that since we do not have the generator commitment history at hand,  we assume a cyclic demand pattern.
To account for this modeling assumption, we define the following index sets for minimum uptime and downtime constraints:
\begin{align*}
& T^{up}_{i,t} := \{x: x= (t-MinUp_{i}+j) \mathrm{mod}\;24 + 1, \\
& \hspace{4em} j \in \mathbb{Z}^{+}_{0},0\leq j \leq MinUp_{i}-1 \},  \forall g \in \mathcal{G}, \forall t \in \mathcal{T}\\
& T^{dw}_{i,t} := \{x: x= (t-MinDw_{i}+j) \mathrm{mod}\;24 + 1, \\
& \hspace{4em} j \in \mathbb{Z}^{+}_{0},0\leq j \leq MinDw_{i}-1 \},  \forall g \in \mathcal{G}, \forall t \in \mathcal{T}.
\end{align*}
Then, we are able to present the formulation as follows:

\begin{subequations}\label{eq:originalminlp}
\begin{align}
\mathrm{min} &\hspace{0.25em}  \sum_{t \in \mathcal{T}} \sum_{i \in \mathcal{G}} (f_{i,t}(u_{i,t},v_{i,t},w_{i,t}) + c_{i,t}(p_{i,t}^g))  \label{eq:UC objective} \\
  \mathrm{s.t.}   &\hspace{0.25em} \text{For each }  i \in \mathcal{G}, t \in \mathcal{T}: \notag \\
  & \hspace{0.25em} u_{i,t'} - u_{i,t} = v_{i,t}- w_{i,t}  \label{eq:logicalequality} \\
  & \hspace{0.25em} v_{i,t} - u_{i,t} \leq 0  \label{eq:logicalineq1} \\
  & \hspace{0.25em} w_{i,t} + u_{i,t} \leq 1 \label{eq:logicalineq2} \\
  & \hspace{0.25em} \sum_{\tau \in T^{up}_{i,t}} v_{i,\tau} \leq u_{i,t} \label{eq:uptime} \\
  & \hspace{0.25em} \sum_{\tau \in T^{dw}_{i,t}} w_{i,\tau} \leq 1- u_{i,t} \label{eq:downtime} \\
  & \hspace{0.25em} p_{i}^{min} u_{i,t} \leq p_{i,t}^{g}\leq  {p}_{i}^{max}u_{i,t}  \label{eq:limitonenergyp} \\
  & \hspace{0.25em} q_{i}^{min} u_{i,t} \leq q_{i,t}^{g}\leq  {q}_{i}^{max}u_{i,t}  \label{eq:limitonenergyq} \\
  & \hspace{0.25em}  u_{i,t},v_{i,t},w_{i,t} \in \{0,1\}  \label{eq:binaries} \\
  & \hspace{0.25em} \eqref{eq:cstheta active balanceu}-\eqref{eq:cstheta reactive balanceu}, 
\eqref{eq:rampopf}-\eqref{eq:vtheta consistency2}.\notag
\end{align}
\end{subequations}
The first component of the objective function \eqref{eq:UC objective} is a linear function of the binary variables: commitment status, startup status and shutdown status of the generators. The second component of the objective function is a function of active power generation of generators, and generally is a quadratic function. Constraints \eqref{eq:logicalequality}, \eqref{eq:logicalineq1} and \eqref{eq:logicalineq2} are logical constraints for the relation between commitment statuses and turn on/off actions of generators. Constraints \eqref{eq:uptime} and \eqref{eq:downtime} enforce minimum uptime and downtime requirements for generator $i$, respectively. Constraints \eqref{eq:limitonenergyp} and \eqref{eq:limitonenergyq} are the limits of real and reactive power outputs of generator $i$, considering the commitment status of the generator. Constraints \eqref{eq:cstheta active balanceu}-\eqref{eq:cstheta reactive balanceu}, 
\eqref{eq:rampopf}-\eqref{eq:vtheta consistency2} account for the AC power flow equations.

\subsubsection {MISOCP Relaxation}
\label{sec:misocp relax}
The MISOCP relaxation of Problem \eqref{eq:originalminlp}  can be obtained 
by relaxing the constraint \eqref{eq:vtheta consistency1} and omitting the constraint \eqref{eq:vtheta consistency2} as previously done in Section \ref{sec:socp relax}.
\begin{subequations}\label{eq:misocprelx}
\begin{align}
\mathrm{min} &\hspace{0.25em}  \sum_{t \in \mathcal{T}} \sum_{i \in \mathcal{G}} (f_{i,t}(u_{i,t},v_{i,t},w_{i,t}) + c_{i,t}(p_{i,t}^g))  \tag{\ref{eq:misocprelx}} \\
  \mathrm{s.t.}& \hspace{0.25em} \eqref{eq:logicalequality}-\eqref{eq:binaries} \notag \\
  &\hspace{0.25em} \eqref{eq:cstheta active balanceu}-\eqref{eq:cstheta reactive balanceu}, 
\eqref{eq:rampopf}-\eqref{eq:ubline bw active},\eqref{eq:conicrelax} \notag.
\end{align}
\end{subequations}

This MISOCP relaxation is heavily used in the solution algorithms that we propose in Section \ref{solutionalgos}.

\subsubsection{DC Approximation}
\label{sec:milp approx}

The formulation for the DC approximation of the UC Problem is given as follows:

\begin{subequations}\label{eq:dcuc}
\begin{align}
\mathrm{min} &\hspace{0.25em}  \sum_{t \in \mathcal{T}} \sum_{i \in \mathcal{G}} (f_{i,t}(u_{i,t},v_{i,t},w_{i,t}) + c_{i,t}(p_{i,t}^g))  \tag{\ref{eq:dcuc}} \\
  \mathrm{s.t.} &\hspace{0.25em} \eqref{eq:logicalequality} - \eqref{eq:limitonenergyp}, \eqref{eq:binaries} \notag \\
  &\hspace{0.25em} \eqref{eq:dcbalance} - \eqref{eq:dcsecurity} \notag.
\end{align}
\end{subequations}

\section{Proposed Solution Algorithms}\label{solutionalgos}

We propose three MISOCP-based algorithms for solving the UC Problem with AC power flow equations: i) base, ii) enhanced and iii) decomposition-based. Base and enhanced algorithms rely on solving the MISOCP relaxations of the UC Problem, and then solving an MOPF Problem to find feasible solutions. The decomposition algorithm is proposed for solving large instances and relies on using Lagrangian decomposition to  solve the time-decomposed MISOCP relaxations of the UC Problem.

\subsection{Base Algorithm}

The base algorithm depends on the observation  that if    the generator commitment schedules are at hand, then the UC Problem reduces to an instance of the  MOPF problem. We utilize this fact and solve an MISOCP relaxation to find a candidate commitment schedule and a lower bound for the UC Problem with AC power follows. Then, we solve  an MOPF Problem with the given commitment status to come up with a feasible solution to the UC Problem. The outline of our approach is given in Algorithm \ref{basealgo}.

\begin{algorithm}[h!]
\label{basealgo}
\SetAlgoLined
\KwResult{Lower and upper bounds for UC problem with AC power flows}
 1. Solve the MISOCP relaxation \eqref{eq:misocprelx} of  Problem \eqref{eq:originalminlp}.\\
 2. Obtain a candidate commitment schedule. \\
 3. Solve Problem \eqref{eq:themopf} with the commitment schedule obtained in the previous step.\\
 4. Calculate the optimality gap.
 \caption{Base Algorithm}
\end{algorithm}

\subsection{Enhanced Algorithm}

In order to improve the lower bound obtained from  Algorithm \ref{basealgo}, we use two different methods taken from \cite{strongsocpkocuk}, which are arctangent envelopes and SDP cut separation. The first method utilizes an outer-approximation of the feasible region of the arctangent constraint by adding linear inequalities. The second method solves an SDP separation problem over some cycles of the network given a solution to the MISOCP relaxation, which is solved by MOSEK \cite{mosek}. From the separation problem, we obtain linear inequalities and add them as a cutting plane to the model. We present the details of this approach in  Algorithm \ref{enhalgo}.

\begin{algorithm}[h!]\label{enhalgo}
\SetAlgoLined
\KwResult{Lower and upper bounds for UC problem with AC power flows}
 1. Compute a cycle basis. \\
 2. For each edge, add the arctangent constraints. \\
 3. Solve the continuous relaxation of  MISOCP \eqref{eq:misocprelx}. \\
 4. For $i=1$ to 5: \\
 \hspace{1.5em}1. Solve the separation problem for each cycle in the cycle basis in parallel. \\
 \hspace{1.5em}2. Add the cuts obtained from the separation problem and resolve the continuous relaxation of  MISOCP \eqref{eq:misocprelx}. \\
 5. Solve the MISOCP problem \eqref{eq:misocprelx} with the cuts obtained. \\
 6. Obtain a candidate commitment schedule. \\
 7. Solve Problem \eqref{eq:themopf} with the commitment schedule obtained in the previous step.\\
 8. Calculate the optimality gap. 
 \caption{Enhanced Algorithm}
\end{algorithm}

\subsection{Decomposition Algorithm}
We develop a decomposition algorithm  since solving the MISOCP problem becomes computationally expensive for large-scale instances. In order to obtain lower bounds faster, we propose dividing the planning horizon into blocks, which are solved separately. For the constraints that include time indices belonging to different blocks (subproblems), we utilize the Lagrangian decomposition method. For deciding on the Lagrangian multipliers, we first solve the continuous relaxation of the MISOCP problem \eqref{eq:misocprelx} and retrieve the optimal dual values of these constraints. We set these dual values of the constraints as the corresponding Lagrangian multipliers and solve the subproblems for each block. In order to ensure that solution of each subproblem outputs a feasible commitment schedule, we solve an intermediate problem, which we call the restricted MISOCP. In this restricted MISOCP, we fix a generator to be turned on/off in the whole planning horizon, if a generator is turned on/off in each of the 24 time periods. Then, for the generators that are turned on for some periods but that are off for some periods, we fix their turned on periods. The solution to this MISOCP becomes our candidate commitment schedule and we solve the MOPF Problem accordingly. We also apply the subgradient algorithm to update the Lagrangian multipliers and try to improve the lower bound.

Before presenting the formulation, we define the following parameters for each generator $i \in \mathcal{G}$, time index $t \in \mathcal{T}$ and block $b$:
\begin{align*}
& 
t'' = 
\begin{cases}
t+1, & \text{if } t \neq 24 \\
1, & \text{if } t = 24
\end{cases}\\
&\underline{t_b} := \min(\mathcal{T}_{b}) \\
&\overline{t_b} := \max(\mathcal{T}_{b}) \\
& T^{intu}_{i,t} := T_{i,t}^{up} \cap (\mathcal{T} \setminus \mathcal{T}_b) \\
& T^{intd}_{i,t} := T_{i,t}^{dw} \cap (\mathcal{T} \setminus \mathcal{T}_b). \\
& T_{i,r1} := \{ t: T_{i,t}^{up} \cap \mathcal{T}_b \neq \emptyset, T_{i,t}^{up} \not\subseteq \mathcal{T}_b \} \\
& T_{i,r2} := \{ t: T^{intu}_{i,t} \neq \emptyset \} \\
& T_{i,r3} := \{ t: T_{i,t}^{dw} \cap \mathcal{T}_b \neq \emptyset, T_{i,t}^{dw} \not\subseteq \mathcal{T}_b \} \\
& T_{i,r4} := \{ t: T^{intd}_{i,t} \neq \emptyset \}.
\end{align*}

In order to formulate the Lagrangian subproblems, we utilize the MISOCP relaxation \eqref{eq:misocprelx} and obtain the subproblems for each block $b$ as follows:

\begin{subequations} \label{eq:lagrrelax}
\begin{align}
\mathrm{min} &\hspace{0.25em}  \sum_{t \in \mathcal{T}_{b}} \sum_{i \in \mathcal{G}} (f_{i,t}(u_{i,t},v_{i,t},w_{i,t}) + c_{i,t}(p_{i,t}^g)) \label{eq:lagrrelaxObj} \\
&\hspace{0.25em} + \sum_{i \in \mathcal{G}} \lambda_{i,\underline{t_b}}^{ru}(p_{i,\underline{t_b}}^{g} - RU_{i}) + \sum_{i \in \mathcal{G}} \lambda_{i,\overline{t_b}''}^{ru}(-p_{i,\overline{t_b}}^{g}) \notag \\
&\hspace{0.25em} + \sum_{i \in \mathcal{G}} \lambda_{i,\underline{t_b}}^{rd}(p_{i,\underline{t_b}}^{g} + RD_{i}) + \sum_{i \in \mathcal{G}} \lambda_{i,\overline{t_b}''}^{rd}(-p_{i,\overline{t_b}}^{g}) \notag \\
&\hspace{0.25em} + \sum _{i \in \mathcal{G}} \lambda_{i,\underline{t_b}}^{log}(-u_{i,\underline{t_b}}+v_{i,\underline{t_b}}-w_{\underline{t_b}}) + \sum _{i \in \mathcal{G}} \lambda_{i,\overline{t_b}''}^{log}(u_{\overline{t_b}}) \notag \\
& \hspace{0.25em}  + \sum_{t \in T_{i,r1}} \sum_{i \in \mathcal{G}} \lambda_{i,t}^{up} ( \sum_{k \in T^{intu}_{i,t}} (v_{i,k}) - u_{i,t}) \notag\\
& \hspace{0.25em}+ \sum_{t \in T_{i,r2}} \sum_{i \in \mathcal{G}} \lambda_{i,t}^{up} ( \sum_{k \in T^{intu}_{i,t}} (v_{i,k}))\notag    \\
& \hspace{0.25em}  + \sum_{t \in T_{i,r3}} \sum_{i \in \mathcal{G}} \lambda_{i,t}^{dw} ( \sum_{k \in T^{intd}_{i,t}} (w_{i,k}) + u_{i,t} -1) \notag\\
& \hspace{0.25em} + \sum_{t \in T_{i,r4}} \sum_{i \in \mathcal{G}} \lambda_{i,t}^{up} ( \sum_{k \in T^{intd}_{i,t}} (w_{i,k}))\notag    \\
  \mathrm{s.t.}   &\hspace{0.25em} \sum_{\tau \in T^{up}_{i,t}} v_{i,\tau} \leq u_{i,t} \hspace{3.5em} i \in \mathcal{G}, t \in \mathcal{T}_b: T^{up}_{i,t} \subseteq \mathcal{T}_b \label{eq:minupu}\\
    &\hspace{0.45em} \sum_{\tau \in T^{dw}_{i,t}} w_{i,\tau} \leq 1-u_{i,t} \hspace{1em} i \in \mathcal{G}, t \in \mathcal{T}_b: T^{dw}_{i,t} \subseteq \mathcal{T}_b \label{eq:mindwu}\\
    &\hspace{0.25em} u_{i,T_{t}} - u_{i,t} + v_{i,t} - w_{i,t} = 0 \hspace{0.75em} i \in \mathcal{G}, t \in \mathcal{T}_b \setminus \{ \underline{t_b}\} \label{eq:logicu} \\
    &\hspace{0.25em} -RD_{i} \leq p^{g}_{i,t} - p^{g}_{i,T_{t}} \hspace{3.6em} i \in \mathcal{G}, t \in \mathcal{T}_b \setminus \{ \underline{t_b}\} \label{eq:rampdwu} \\
    &\hspace{0.25em} p^{g}_{i,t} - p^{g}_{i,T_{t}} \leq RU_{i} \hspace{5em} i \in \mathcal{G}, t \in \mathcal{T}_b \setminus \{ \underline{t_b}\} \label{eq:rampupu}\\
   & \eqref{eq:limitonenergyp}-\eqref{eq:binaries}, \eqref{eq:cstheta active balanceu}-\eqref{eq:cstheta reactive balanceu}, 
\eqref{eq:rampopf}-\eqref{eq:ubline bw active},\eqref{eq:conicrelax}.\notag
\end{align}
\end{subequations}

Our decomposition method is outlined in Algorithm \ref{decompalgo}.
\begin{algorithm}[h!]\label{decompalgo}
\SetAlgoLined
\KwResult{Lower and upper bounds for UC with AC power flows}
 1. Compute a cycle basis. \\
 2. For each edge, add the arctangent constraints. \\
 3. Solve the continuous relaxation of the MISOCP \eqref{eq:misocprelx}. \\
 4. For $i=1 $ to 5: \\
 \hspace{1.5em}1. Solve the separation problem for each cycle in the cycle basis in parallel. \\
 \hspace{1.5em}2. Add the cuts obtained from the separation problem and resolve the continuous relaxation of MISOCP \eqref{eq:misocprelx}. \\
 5. Obtain the Lagrangian multipliers from the previous step. \\
 6. For $i=1$ to 5\\
 \hspace{1.5em}1. For each block $b$, solve the Lagrangian Relaxation Problem \eqref{eq:lagrrelax}.\\
 \hspace{1.5em}2. Obtain the commitment decisions for the original time horizon, solve the restricted MISOCP.\\
 \hspace{1.5em}3. Solve Problem \eqref{eq:themopf} in order to find a feasible solution to the UC with AC power flows problem.\\
 \hspace{1.5em}4. Update the Lagrangian Multiplier values with the subgradient algorithm.\\
 7. Take the best lower bound and the best feasible solution as the upper bound from the previous step.\\
 8. Calculate the optimality gap.
 \caption{Decomposition Algorithm}
\end{algorithm}

\section{Computational Results}\label{compresults}

\subsection{Instance Creation}
We face some  difficulties in finding realistic problem instances for the UC Problem with AC power flows. Therefore, we decide to create our own problem instances based on NESTA AC OPF instances \cite{nestaref}. To be able to convert the AC OPF instance to an AC UC instance, we need the following parameters for each generator $i \in \mathcal{G}$: ramp-up/down rate, minimum up/down time, startup cost and fixed cost. Also, for each bus $i \in \mathcal{B}$, we need the demand for 24 periods. Due to lack of generator history, we assumed that the demand is cyclic. The details of the creation procedure are provided in  Appendix \ref{app:instance}.  Instances created according to this procedure can be found at the following website: \url{https://sites.google.com/site/burakkocuk/research}.


\subsection{Computational Setting}

The experiments are carried out on a desktop workstation with 3.7 GHz processor and 32 GB of RAM. Lower bounds are obtained by solving MISOCP, utilizing Gurobi \cite{gurobi}. Feasible solutions are obtained by IPOPT \cite{ipoptrefer}, a local solver, which provides upper bounds. We solve three types of instances which we denote by TYP, API and SAD and they refer to typical, congested and small angle difference condition instances, respectively. \texttt{DC}, \texttt{MISOCP} and \texttt{MISOCP++} refer to the DC approach, base algorithm and enhanced algorithm, respectively. LBT and UBT refer to the time to find lower bound and upper bound, respectively. LB and UB refers to the lower bound and upper bound, respectively. `local inf' denotes that IPOPT converged to a locally infeasible point. We calculate \%Gap  according to the following formula: $\frac{UB - LB}{UB} \times 100$. We set the relative  optimality tolerance of Gurobi as 0.1\% unless otherwise stated.

\subsection{DC vs MISOCP Based Methods}

We  first run some preliminary experiments to compare the DC and MISOCP based methods. 
 The computational results presented in Table \ref{tab:dcvsmisocp} show that the \texttt{DC} approach generally finds a solution in a shorter time than the MISOCP based methods (with some abuse of terminology, we report the computational time of solving the DC based method under the LBT column). However, the solution of Problem \eqref{eq:dcapprox} might not be feasible, it may not lead to an AC-feasible generator commitment schedule  or it may lead to a suboptimal AC-feasible  generator commitment schedule as the results obtained from the reported three instances respectively show. Therefore, we only focus on the MISOCP based solution methods in the remainder of the paper since they are consistently more successful.


\begin{table}[h!]\caption{DC vs. MISOCP} \label{tab:dcvsmisocp}
\centering
\begin{tabular}{|llrrrr|}
\hline
Case & Method & LBT (s) & UBT (s)& UB &\%Gap \\ \hline
case9-SAD & \texttt{DC}       & infeas. & -     & -        & -    \\
case9-SAD& \texttt{MISOCP}   & 0.16    & 0.77  & 1582.36  & 0.52    \\
case9-SAD & \texttt{MISOCP++} & 2.24   & 0.76  & 1582.36   & 0.03 \\ \hline
case14-TYP& \texttt{DC}       & 0.01    & -  & infeas.  & -    \\
case14-TYP & \texttt{MISOCP}   & 1.57   & 1.67 & 227.14  & 0.10    \\
case14-TYP & \texttt{MISOCP++} & 8.93  & 1.24  & 227.13 & 0.01 \\ \hline
case57-SAD & \texttt{DC}       & 0.11    & 6.39  & 642.30  & -    \\
case57-SAD & \texttt{MISOCP}   & 7.08   & 6.21  & 545.44   & 0.67 \\
case57-SAD & \texttt{MISOCP++} & 33.88  & 6.43  & 545.44   & 0.19 \\ \hline
\end{tabular}
\end{table}

\subsection{Computational Results for the Base and Enhanced Algorithms}

In Table \ref{tab:enhancedresults}, we present the computational results for the base and the enhanced algorithm.

\begin{table*}[]\caption{Results for \texttt{MISOCP} and \texttt{MISOCP++}} \label{tab:enhancedresults}
\centering
\begin{tabular}{|l|l|rrrr|rrrr|}
\hline
              &                           & \multicolumn{4}{c|}{\texttt{MISOCP}}           & \multicolumn{4}{c|}{\texttt{MISOCP++}}         \\ \hline
instance      & \multicolumn{1}{r|}{Type} & LBT (s)  & UBT (s)  & UB        & \%Gap & LBT (s)  & UBT (s)  & UB        & \%Gap \\ \hline\hline
6ww       & TYP                       & 0.17    & 0.74    & 4417.38   & 0.02  & 5.02    & 0.73    & 4417.38   & 0.00  \\
9wscc     & TYP                       & 0.16    & 0.74    & 1563.87   & 0.00  & 2.23    & 0.68    & 1563.87   & 0.00  \\
14ieee    & TYP                       & 1.57    & 1.67    & 227.14    & 0.10  & 8.93    & 1.24    & 227.13    & 0.01  \\
24ieeerts & TYP                       & 0.96    & 3.56    & 155161.46 & 0.00  & 14.24   & 3.51    & 155161.46 & 0.00  \\
30as      & TYP                       & 0.84    & 2.32    & 1832.06   & 0.00  & 13.81   & 2.30    & 1832.06   & 0.00  \\
30ieee    & TYP                       & 6.87    & 3.14    & 235.46    & 2.14  & 20.66   & 2.92    & 235.45    & 0.02  \\
39epri    & TYP                       & 39.82   & 4.57    & 28336.21  & 0.03  & 73.21   & 4.55    & 28336.21  & 0.12  \\
57ieee    & TYP                       & 14.45   & 5.84    & 542.08    & 0.08  & 39.59   & 5.95    & 542.08    & 0.06  \\ \hline\hline
6ww       & API                       & 0.20    & 0.89    & 421.97    & 0.21  & 5.24    & 0.84    & 421.97    & 0.05  \\
9wscc     & API                       & 0.17    & 0.83    & 434.30    & 0.08  & 2.17    & 0.80    & 434.30    & 0.00  \\
14ieee    & API                       & 2.12    & 1.53    & 214.07    & 0.13  & 8.45    & 1.36    & 214.06    & 0.02  \\
24ieeerts & API                       & 0.93    & 3.57    & 11335.96  & 1.23  & 15.07   & 3.55    & 11335.96  & 0.74  \\
30as      & API                       & 0.90    & 2.34    & 795.26    & 0.31  & 15.73   & 2.37    & 795.26    & 0.18  \\
30ieee    & API                       & 6.58    & 4.79    & local inf & N/A   & 31.36   & 3.37    & 273.31    & 0.00  \\
39epri    & API                       & 18.63   & 4.73    & 2211.70   & 1.34  & 49.16   & 4.73    & 2211.70   & 1.30  \\
57ieee    & API                       & 17.35   & 10.92   & 632.30    & 0.05  & 49.77   & 11.04   & 632.30    & 0.02  \\ \hline\hline
6ww       & SAD                       & 0.16    & 0.67    & 4417.40   & 0.02  & 5.07    & 0.72    & 4417.40   & 0.00  \\
9wscc     & SAD                       & 0.16    & 0.77    & 1582.36   & 0.52  & 2.24    & 0.76    & 1582.36   & 0.03  \\
14ieee    & SAD                       & 1.08    & 1.51    & local inf & N/A   & 8.48    & 1.48    & 227.13    & 0.01  \\
24ieeerts & SAD                       & 3.02    & 4.84    & 157448.37 & 0.71  & 18.56   & 4.80    & 157448.37 & 0.25  \\
30as      & SAD                       & 0.77    & 2.60    & 1843.64   & 0.48  & 16.68   & 2.56    & 1843.64   & 0.02  \\
30ieee    & SAD                       & 5.20    & 2.91    & 235.44    & 1.02  & 20.59   & 3.11    & 235.44    & 0.01  \\
39epri    & SAD                       & 36.50   & 5.89    & local inf & N/A   & 68.63   & 3.38    & 29687.04  & 0.13  \\
57ieee    & SAD                       & 7.08    & 6.21    & 545.44    & 0.67  & 33.88   & 6.43    & 545.44    & 0.19  \\ \hline
\end{tabular}
\end{table*}

The \texttt{MISOCP} method fails to produce a feasible commitment schedule in one of the congested instances and two of the small angle difference instances, whereas \texttt{MISOCP++} method is able to provide a feasible commitment schedule for each instance. Out of the instances that a feasible commitment schedule is obtained, \texttt{MISOCP} method and \texttt{MISOCP++} method have 0.43\% and 0.13\% optimality gap in average, respectively.
For TYP instances, the \texttt{MISOCP++} outperforms \texttt{MISOCP}, by reducing average optimality gap from 0.29\% to 0.03\%.   For API instances, the \texttt{MISOCP} method is able to solve seven out of eight instances with the average optimality gap of 0.47\% whereas the \texttt{MISOCP++} method yields a feasible solution to all eight instances with average optimality gap of 0.29\%. Especially for the SAD instances, six out of eight instances could be solved with \texttt{MISOCP} with an average optimality gap of 0.57\%, whereas \texttt{MISOCP++} method yields feasible solutions to all eight instances with an average optimality gap of 0.08\%.

The average time to obtain a lower bound for the problem using \texttt{MISOCP} and \texttt{MISOCP++} are 6.9 and 22 seconds, respectively. An increase in the time to obtain a lower bound is expected since we are solving SDP separation problems in addition to the SOCP relaxations in the latter method. However, the additional computational effort is justified since we obtain stronger lower bounds and, hence, smaller optimality gaps. 

We also experiment with different cost combinations than reported in the paper and obtain similarly successful results \cite{denizthesis}.

We see that all the problem instances were solved to a less than 1.3\% optimality gap and within a reasonable timespan, considering the challenging nature of the problem. However, for larger instances, finding a lower bound is time consuming. Therefore, we came up with a decomposition method in order to solve larger instances.

\subsection{Decomposition Method}

In Table \ref{tab:decompositionresults}, we present the computational results for the decomposition algorithm. Lagrangian subproblems are solved to 1\% optimality tolerances. We note that we terminate the subgradient iterations once a feasible solution of the UC Problem with 1\% proven optimality gap is obtained.

\begin{table*}[h!]\caption{Results for the Decomposition Algorithm} \label{tab:decompositionresults}
\centering
\begin{tabular}{|lc|r|rrrrr|rrrr|} 
\hline
Cases                    & Type                 & $b$ & LBT (s)  & UBT (s) & LB        & UB        & \%Gap & Total Time (s) & LB       & UB       & \%Gap  \\ 
\hline
                         &                      &   & \multicolumn{5}{c|}{After 1 iteration of the subgradient algorithm}           & \multicolumn{4}{c|}{After 5 iterations of the subgradient algorithm}      \\ 
\hline
\hline
\multirow{2}{*}{57ieee}  & \multirow{2}{*}{TYP} & 4 & 36.12   & 7.39   & 541.94    & 542.08    & 0.03  & ~             & ~        & ~        & ~      \\
                         &                      & 6 & 33.92   & 7.38   & 541.94    & 542.08    & 0.03  & ~             & ~        & ~        & ~      \\ 
\hline
\multirow{2}{*}{89peg}   & \multirow{2}{*}{TYP} & 4 & 21.74   & 98.17  & 2343.29   & 2344.56   & 0.05  & ~             & ~        & ~        & ~      \\
                         &                      & 6 & 24.12   & 102.41 & 2343.29   & 2344.56   & 0.05  & ~             & ~        & ~        & ~      \\ 
\hline
\multirow{2}{*}{118ieee} & \multirow{2}{*}{TYP} & 4 & 1351.24 & 24.13  & 1550.08   & 1878.93   & 17.50 & 4958.19       & 1551.65  & 1722.95  & 9.94   \\
                         &                      & 6 & 542.61  & 34.41  & 1558.35   & 2232.58   & 30.20 & 2057.51       & 1558.35  & 1961.92  & 20.57  \\ 
\hline
\multirow{2}{*}{500goc}  & \multirow{2}{*}{TYP} & 4 & 744.13  & 384.71 & 626725.03 & 626791.30 & 0.01  & ~             & ~        & ~        & ~      \\
                         &                      & 6 & 717.93  & 387.07 & 626501.81 & 626791.30 & 0.05  & ~             & ~        & ~        & ~      \\ 
\hline
\hline
\multirow{2}{*}{57ieee}  & \multirow{2}{*}{API} & 4 & 35.17   & 5.30   & 631.61    & 632.61    & 0.16  & ~             & ~        & ~        & ~      \\
                         &                      & 6 & 33.43   & 5.18   & 631.98    & 632.61    & 0.10  & ~             & ~        & ~        & ~      \\ 
\hline
\multirow{2}{*}{89peg}   & \multirow{2}{*}{API} & 4 & 117.13  & 325.12 & 49896.65  & 52984.97  & 5.83  & 2075.29       & 49931.52 & 52984.97 & 5.76   \\
                         &                      & 6 & 35.70   & 231.44 & 49100.91  & 54954.02  & 10.65 & 1309.04       & 49465.83 & 54954.02 & 9.99   \\ 
\hline
\multirow{2}{*}{500goc}  & \multirow{2}{*}{API} & 4 & 548.13  & 154.03 & 663900.96 & 664777.50 & 0.13  & ~             & ~        & ~        & ~      \\
                         &                      & 6 & 547.41  & 153.83 & 663902.33 & 664777.50 & 0.13  & ~             & ~        & ~        & ~      \\ 
\hline
\hline
\multirow{2}{*}{57ieee}  & \multirow{2}{*}{SAD} & 4 & 29.44   & 5.30   & 543.85    & 545.64    & 0.33  & ~             & ~        & ~        & ~      \\
                         &                      & 6 & 26.61   & 5.43   & 543.88    & 545.64    & 0.32  & ~             & ~        & ~        & ~      \\ 
\hline
\multirow{2}{*}{89peg}   & \multirow{2}{*}{SAD} & 4 & 22.50   & 87.12  & 53567.31  & 53584.94  & 0.03  & ~             & ~        & ~        & ~      \\
                         &                      & 6 & 21.84   & 84.71  & 53567.31  & 53584.94  & 0.03  & ~             & ~        & ~        & ~      \\ 
\hline
\multirow{2}{*}{118ieee} & \multirow{2}{*}{SAD} & 4 & 621.12  & 36.75  & 1786.13   & 2292.07   & 22.07 & 3737.20       & 1786.13  & 2292.07  & 22.07  \\
                         &                      & 6 & 346.83  & 20.14  & 1738.59   & 2420.72   & 28.18 & 1613.36       & 1743.32  & 2341.15  & 25.54  \\ 
\hline
\multirow{2}{*}{500goc}  & \multirow{2}{*}{SAD} & 4 & 569.36  & 833.71 & 626560.53 & 627438.64 & 0.14  & ~             & ~        & ~        & ~      \\
                         &                      & 6 & 601.32  & 794.42 & 626549.95 & 627438.64 & 0.14  & ~             & ~        & ~        & ~      \\
\hline
\end{tabular}
\end{table*}

The decomposition method yields promising results for larger instances. 
For the instances 57ieee, 89peg and 500goc, we are able to find feasible solutions with less than 1\% optimality gap (except for case 89peg-API, in which the optimality gap is more than 5\%). 

For the 118-bus instances, which are the most challenging instances we consider, our decomposition method manages to output a feasible solution. For these instances, the base and enhanced algorithms fail to produce a feasible solution with the commitment schedule we obtain within an hour. For the case 118ieee-TYP, the subgradient iterations lead to an improved feasible solution and reduced the optimality gap from 17.50\% to 9.94\%. For the case 118ieee-SAD, we have a relatively large optimality gap, due to the operating conditions and density of the network. We would like to note that such larger optimality gaps are not rare in the OPF literature \cite{babaeinejadsarookolaee2019power}.

In general, when we divide the planning horizon into four blocks, we achieve better optimality gaps compared to dividing the planning horizon into six blocks, with the burden of spending more time to achieve a feasible solution. In addition, applying the decomposition method to the 57-bus instances  results in solving the problem in less time, while achieving similarly successful objective values compared to the \texttt{MISOCP++} method.

Finally, we note that we face numerical difficulties  for the  118ieee-API instance and are not able to obtain a feasible solution.  Also, for 89peg instances, we only add arctangent envelopes but not SDP inequalities, due to numerical difficulties.

\section{Conclusion}\label{conclusion}
In this paper, we  studied the UC Problem with AC power flow equations. We developed three solution methods that utilize the MISOCP relaxation of this challenging problem. Since there was a lack of publicly available problem instances, we  constructed realistic  UC Problem instances with AC power flows  based on publicly available AC OPF instances.
We empirically showed that the classical DC-based approach might yield unfavorable generator commitment schedules whereas our MISOCP-based solution approaches consistently provided provably high-quality feasible solutions over these instances. In particular,  we were able to provide a feasible solution within 1.3\%  optimality gap over instances with up to 57-bus and 24-hour planning horizon in at most two minutes of computational time.
To solve even larger instances, we  developed a decomposition method based on the Lagrangian relaxation approach. The decomposition method enabled us to solve instances with up to 500 buses. 

Future work may include strengthening the decomposition method and solving instances with more than 500 buses. For such larger scale instances, we anticipate that both temporal and spatial decompositions are needed.


%

\appendices
\section{Instance Creation Details}
\label{app:instance}

For each bus, we utilize the real demand profiles taken from \cite{profile1}, \cite{profile2} and \cite{anyacastillo} and randomly assign these profiles. We assume that peak demand of the demand profile is equal to the demand of the AC OPF instance and normalize the demand accordingly. For the reactive demand, we utilize the reactive demand from \cite{anyacastillo}. The demand profiles are given in the Table \ref{tab:demprof}. 

\begin{table}[h!]
\caption{Demand profiles for creation of the problem instances}\label{tab:demprof}
\centering
\begin{tabular}{|c|ccccc|}
\hline
\textbf{\begin{tabular}[c]{@{}c@{}}Time \\ period\end{tabular}} &
  \textbf{\begin{tabular}[c]{@{}c@{}}Real\\ Profile 1\end{tabular}} &
  \textbf{\begin{tabular}[c]{@{}c@{}}Real\\ Profile 2\end{tabular}} &
  \textbf{\begin{tabular}[c]{@{}c@{}}Real\\ Profile 3\end{tabular}} &
  \textbf{\begin{tabular}[c]{@{}c@{}}Max Real\\ Profile\end{tabular}} &
  \textbf{\begin{tabular}[c]{@{}c@{}}Reactive\\ Profile\end{tabular}} \\ \hline
\textbf{1}  & 0.68 & 0.57 & 0.67 & 0.68 & 0.68 \\
\textbf{2}  & 0.64 & 0.64 & 0.63 & 0.64 & 0.65 \\
\textbf{3}  & 0.61 & 0.68 & 0.60 & 0.68 & 0.62 \\
\textbf{4}  & 0.60 & 0.71 & 0.59 & 0.71 & 0.60 \\
\textbf{5}  & 0.60 & 0.75 & 0.59 & 0.75 & 0.61 \\
\textbf{6}  & 0.62 & 0.78 & 0.60 & 0.78 & 0.63 \\
\textbf{7}  & 0.67 & 0.82 & 0.74 & 0.82 & 0.68 \\
\textbf{8}  & 0.74 & 0.85 & 0.86 & 0.86 & 0.69 \\
\textbf{9}  & 0.80 & 0.88 & 0.95 & 0.95 & 0.73 \\
\textbf{10} & 0.84 & 0.92 & 0.96 & 0.96 & 0.81 \\
\textbf{11} & 0.89 & 0.97 & 0.96 & 0.97 & 0.89 \\
\textbf{12} & 0.92 & 1.00 & 0.95 & 1.00 & 0.92 \\
\textbf{13} & 0.94 & 0.92 & 0.95 & 0.95 & 0.95 \\
\textbf{14} & 0.95 & 0.88 & 0.95 & 0.95 & 0.95 \\
\textbf{15} & 0.97 & 0.85 & 0.93 & 0.97 & 0.97 \\
\textbf{16} & 0.99 & 0.78 & 0.94 & 0.99 & 1.00 \\
\textbf{17} & 1.00 & 0.71 & 0.99 & 1.00 & 1.00 \\
\textbf{18} & 0.96 & 0.78 & 1.00 & 1.00 & 0.96 \\
\textbf{19} & 0.96 & 0.85 & 1.00 & 1.00 & 0.96 \\
\textbf{20} & 0.92 & 0.92 & 0.96 & 0.96 & 0.93 \\
\textbf{21} & 0.92 & 0.85 & 0.91 & 0.92 & 0.93 \\
\textbf{22} & 0.88 & 0.78 & 0.83 & 0.88 & 0.91 \\
\textbf{23} & 0.78 & 0.71 & 0.73 & 0.78 & 0.77 \\
\textbf{24} & 0.76 & 0.64 & 0.63 & 0.76 & 0.76 \\ \hline
\end{tabular}
\end{table}

For the determination of UC-specific parameters for each generator, we randomly assign each generator a type, and then calculate the parameters based on their types. The calculation of the parameters are described as follows:
For each generator $i \in \mathcal{G}$, we assign three types, Type 1, 2 and 3.\\
If generator $i$ is of Type 1, then $RU_{i}=RD_{i}=\max \{ p^{min}_i,\frac{p^{max}_{i}}{2}\}$, $MinUp_{i} = MinDw_{i} = 2$,\\
If generator $i$ is of Type 2, then $RU_{i}=RD_{i}=\max \{ p^{min}_i,\frac{p^{max}_{i}}{3}\}$, $MinUp_{i} = MinDw_{i} = 3$,\\
If generator $i$ is of Type 3, then $RU_{i}=RD_{i}=\max \{ p^{min}_i,\frac{p^{max}_{i}}{5}\}$, $MinUp_{i} = MinDw_{i} = 4$.

We choose the following cost combinations for the fixed and startup costs of the generators, respectively:  $F_{i} = 5Li_{i}$ and $StUp_{i} =100Li_{i}$, where $Li_{i}$ is the linear cost of power generation. The shutdown cost is selected as zero in accordance with the literature.

\section*{Acknowledgment}

The authors would like to thank The Scientific and Technological Research Council of Turkey (TÜBİTAK) for supporting this study with project number 119M855.

\ifCLASSOPTIONcaptionsoff
  \newpage
\fi



\bibliographystyle{IEEEtran}
\bibliography{MA_Library.bib}

\begin{thebibliography}{10}
\providecommand{\url}[1]{#1}
\csname url@samestyle\endcsname
\providecommand{\newblock}{\relax}
\providecommand{\bibinfo}[2]{#2}
\providecommand{\BIBentrySTDinterwordspacing}{\spaceskip=0pt\relax}
\providecommand{\BIBentryALTinterwordstretchfactor}{4}
\providecommand{\BIBentryALTinterwordspacing}{\spaceskip=\fontdimen2\font plus
\BIBentryALTinterwordstretchfactor\fontdimen3\font minus
  \fontdimen4\font\relax}
\providecommand{\BIBforeignlanguage}[2]{{%
\expandafter\ifx\csname l@#1\endcsname\relax
\typeout{** WARNING: IEEEtran.bst: No hyphenation pattern has been}%
\typeout{** loaded for the language `#1'. Using the pattern for}%
\typeout{** the default language instead.}%
\else
\language=\csname l@#1\endcsname
\fi
#2}}
\providecommand{\BIBdecl}{\relax}
\BIBdecl

\bibitem{frank1}
S.~Frank, I.~Steponavičė, and S.~Rebennack, ``Optimal power flow: a
  bibliographic survey i,'' \emph{Energy Systems}, vol.~3, 09 2012.

\bibitem{sasson}
A.~M. Sasson, F.~Viloria, and F.~Aboytes, ``Optimal load flow solution using
  the hessian matrix,'' \emph{IEEE Transactions on Power Apparatus and
  Systems}, vol. PAS-92, no.~1, pp. 31--41, 1973.

\bibitem{dacosta}
G.~da~Costa, ``Optimal reactive dispatch through primal-dual method,''
  \emph{IEEE Transactions on Power Systems}, vol.~12, no.~2, pp. 669--674,
  1997.

\bibitem{tognola}
G.~Tognola and R.~Bacher, ``Unlimited point algorithm for opf problems,''
  \emph{IEEE Transactions on Power Systems}, vol.~14, no.~3, pp. 1046--1054,
  1999.

\bibitem{waqquas}
W.~Bukhsh, A.~Grothey, K.~McKinnon, and P.~Trodden, ``Local solutions of the
  optimal power flow problem,'' \emph{Power Systems, IEEE Transactions on},
  vol.~28, pp. 4780--4788, 11 2013.

\bibitem{kocuk016}
\BIBentryALTinterwordspacing
B.~Kocuk, S.~S. Dey, and X.~A. Sun, ``Inexactness of sdp relaxation and valid
  inequalities for optimal power flow,'' \emph{IEEE Transactions on Power
  Systems}, vol.~31, no.~1, p. 642–651, Jan 2016. [Online]. Available:
  \url{http://dx.doi.org/10.1109/TPWRS.2015.2402640}
\BIBentrySTDinterwordspacing

\bibitem{lpapproxcoffrin}
C.~Coffrin and P.~Van~Hentenryck, ``A linear-programming approximation of ac
  power flows,'' \emph{INFORMS Journal on Computing}, vol.~26, 06 2012.

\bibitem{stottdc}
B.~Stott, J.~Jardim, and O.~Alsac, ``Dc power flow revisited,'' \emph{Power
  Systems, IEEE Transactions on}, vol.~24, pp. 1290 -- 1300, 09 2009.

\bibitem{bienstock}
D.~Bienstock and G.~Munoz, ``On linear relaxations of opf problems,'' 2014.

\bibitem{griffith}
\BIBentryALTinterwordspacing
R.~E. Griffith and R.~A. Stewart, ``A nonlinear programming technique for the
  optimization of continuous processing systems,'' \emph{Management Science},
  vol.~7, no.~4, pp. 379--392, 1961. [Online]. Available:
  \url{https://doi.org/10.1287/mnsc.7.4.379}
\BIBentrySTDinterwordspacing

\bibitem{baiwei1}
X.~Bai and H.~Wei, ``Semi-definite programming-based method for
  security-constrained unit commitment with operational and optimal power flow
  constraints,'' \emph{Generation, Transmission and Distribution, IET}, vol.~3,
  pp. 182 -- 197, 03 2009.

\bibitem{baiwei2}
X.~Bai, H.~Weihua, K.~Fujisawa, and Y.~Wang, ``Semidefinite programming for
  optimal power flow problems,'' \emph{International Journal of Electrical
  Power and Energy Systems}, vol.~30, pp. 383--392, 07 2008.

\bibitem{lavaei}
J.~Lavaei and S.~H. Low, ``Zero duality gap in optimal power flow problem,''
  \emph{IEEE Transactions on Power Systems}, vol.~27, no.~1, pp. 92--107, 2012.

\bibitem{madani}
R.~Madani, M.~Ashraphijuo, and J.~Lavaei, ``Promises of conic relaxation for
  contingency-constrained optimal power flow problem,'' \emph{IEEE Transactions
  on Power Systems}, vol.~31, no.~2, pp. 1297--1307, 2016.

\bibitem{jabrcone}
R.~Jabr, ``Radial distribution load flow using conic programming,'' \emph{IEEE
  Transactions on Power Systems}, vol.~21, no.~3, pp. 1458--1459, 2006.

\bibitem{strongsocpkocuk}
B.~Kocuk, S.~Dey, and X.~Sun, ``Strong socp relaxations for the optimal power
  flow problem,'' \emph{Operations Research}, vol.~64, 05 2016.

\bibitem{alguacil}
N.~Alguacil and A.~Conejo, ``Multiperiod optimal power flow using benders
  decomposition,'' \emph{IEEE Transactions on Power Systems}, vol.~15, no.~1,
  pp. 196--201, 2000.

\bibitem{demirovic}
N.~Demirovic, S.~Tesnjak, and A.~Tokic, ``Hot start and warm start in lp based
  interior point method and it's application to multiperiod optimal power
  flows,'' in \emph{2006 IEEE PES Power Systems Conference and Exposition},
  2006, pp. 699--704.

\bibitem{schanen}
M.~Schanen, F.~Gilbert, C.~G. Petra, and M.~Anitescu, ``Toward multiperiod
  ac-based contingency constrained optimal power flow at large scale,'' in
  \emph{2018 Power Systems Computation Conference (PSCC)}, 2018, pp. 1--7.

\bibitem{muckstadt}
\BIBentryALTinterwordspacing
J.~A. Muckstadt and S.~A. Koenig, ``An application of lagrangian relaxation to
  scheduling in power-generation systems,'' \emph{Operations Research},
  vol.~25, no.~3, pp. 387--403, 1977. [Online]. Available:
  \url{https://doi.org/10.1287/opre.25.3.387}
\BIBentrySTDinterwordspacing

\bibitem{dubost}
\BIBentryALTinterwordspacing
L.~Dubost, R.~Gonzalez, and C.~Lemar{\'{e}}chal, ``A primal-proximal heuristic
  applied to the french unit-commitment problem,'' \emph{Math. Program.}, vol.
  104, no.~1, pp. 129--151, 2005. [Online]. Available:
  \url{https://doi.org/10.1007/s10107-005-0593-4}
\BIBentrySTDinterwordspacing

\bibitem{beltran02}
C.~Beltran and F.-J. Heredia, ``Unit commitment by augmented lagrangian
  relaxation: testing two decomposition approaches,'' \emph{Journal of
  optimization theory and applications}, vol. 112, no.~2, pp. 295--314, Feb
  2002.

\bibitem{carroyo}
M.~Carrion and J.~Arroyo, ``A computationally efficient mixed-integer linear
  formulation for the thermal unit commitment problem,'' \emph{IEEE
  Transactions on Power Systems}, vol.~21, no.~3, pp. 1371--1378, 2006.

\bibitem{morales13}
G.~Morales-España, J.~M. Latorre, and A.~Ramos, ``Tight and compact milp
  formulation of start-up and shut-down ramping in unit commitment,'' in
  \emph{2013 IEEE Power Energy Society General Meeting}, 2013, pp. 1--1.

\bibitem{frangioni09}
F.~Antonio, C.~Gentile, and F.~Lacalandra, ``Tighter approximated milp
  formulations for unit commitment problems,'' \emph{Power Systems, IEEE
  Transactions on}, vol.~24, pp. 105 -- 113, 03 2009.

\bibitem{ostrowski12}
J.~Ostrowski, M.~Anjos, and A.~Vannelli, ``Tight mixed integer linear
  programming formulations for the unit commitment problem,'' \emph{IEEE
  Transactions on Power Systems - IEEE TRANS POWER SYST}, vol.~27, pp. 39--46,
  02 2012.

\bibitem{nanou21}
\BIBentryALTinterwordspacing
S.~I. Nanou, G.~N. Psarros, and S.~A. Papathanassiou, ``Network-constrained
  unit commitment with piecewise linear ac power flow constraints,''
  \emph{Electric Power Systems Research}, vol. 195, p. 107125, 2021. [Online].
  Available:
  \url{https://www.sciencedirect.com/science/article/pii/S0378779621001061}
\BIBentrySTDinterwordspacing

\bibitem{fuetal}
Y.~Fu, M.~Shahidehpour, and Z.~Li, ``Ac contingency dispatch based on
  security-constrained unit commitment,'' \emph{IEEE Transactions on Power
  Systems}, vol.~21, no.~2, pp. 897--908, 2006.

\bibitem{sifuentes}
W.~Sifuentes and A.~Vargas, ``Hydrothermal scheduling using benders
  decomposition: Accelerating techniques,'' \emph{Power Systems, IEEE
  Transactions on}, vol.~22, pp. 1351 -- 1359, 09 2007.

\bibitem{bendersslow}
\BIBentryALTinterwordspacing
D.~W. Watkins and D.~C. McKinney, ``Decomposition methods for water resources
  optimization models with fixed costs,'' \emph{Advances in Water Resources},
  vol.~21, no.~4, pp. 283--295, 1998. [Online]. Available:
  \url{https://www.sciencedirect.com/science/article/pii/S0309170896000619}
\BIBentrySTDinterwordspacing

\bibitem{atamturk}
S.~Fattahi, M.~Ashraphijuo, J.~Lavaei, and A.~Atamturk, ``Conic relaxation of
  the unit commitment problem,'' \emph{Energy}, vol. 134, 10 2016.

\bibitem{paredes}
\BIBentryALTinterwordspacing
M.~Paredes, L.~Martins, S.~Soares, and H.~Ye, ``Benders’ decomposition of the
  unit commitment problem with semidefinite relaxation of ac power flow
  constraints,'' \emph{Electric Power Systems Research}, vol. 192, p. 106965,
  Mar 2021. [Online]. Available:
  \url{http://dx.doi.org/10.1016/j.epsr.2020.106965}
\BIBentrySTDinterwordspacing

\bibitem{anyacastillo}
A.~Castillo, C.~Laird, C.~A. Silva-Monroy, J.-P. Watson, and R.~P. O’Neill,
  ``The unit commitment problem with ac optimal power flow constraints,''
  \emph{IEEE Transactions on Power Systems}, vol.~31, no.~6, pp. 4853--4866,
  2016.

\bibitem{liuetal}
J.~Liu, C.~D. Laird, J.~K. Scott, J.~P. Watson, and A.~Castillo, ``Global
  solution strategies for the network-constrained unit commitment problem with
  ac transmission constraints,'' \emph{IEEE Transactions on Power Systems},
  vol.~34, no.~2, 10 2018.

\bibitem{exposito}
A.~Gómez~Expósito and E.~Romero~Ramos, ``Reliable load flow technique for
  radial distribution networks,'' \emph{IEEE Transactions on Power Systems},
  vol.~14, no.~3, pp. 1063--1069, 1999.

\bibitem{mosek}
\BIBentryALTinterwordspacing
M.~ApS, \emph{MOSEK Optimizer API for Python}, 2021. [Online]. Available:
  \url{https://docs.mosek.com/8.1/pythonapi/index.html}
\BIBentrySTDinterwordspacing

\bibitem{nestaref}
C.~Coffrin, D.~Gordon, and P.~Scott, ``Nesta, the nicta energy system test case
  archive,'' 2019.

\bibitem{gurobi}
\BIBentryALTinterwordspacing
L.~Gurobi~Optimization, ``Gurobi optimizer reference manual,'' 2021. [Online].
  Available: \url{http://www.gurobi.com}
\BIBentrySTDinterwordspacing

\bibitem{ipoptrefer}
A.~Wächter and L.~Biegler, ``On the implementation of an interior-point filter
  line-search algorithm for large-scale nonlinear programming,''
  \emph{Mathematical programming}, vol. 106, pp. 25--57, 03 2006.

\bibitem{denizthesis}
D.~Tuncer, ``Misocp-based solution approaches to the unit commitment problem
  with ac power flows,'' Master's thesis, Sabanci University, Industrial
  Engineering, Istanbul, 2021.

\bibitem{babaeinejadsarookolaee2019power}
S.~Babaeinejadsarookolaee, A.~Birchfield, R.~D. Christie, C.~Coffrin,
  C.~DeMarco, R.~Diao, M.~Ferris, S.~Fliscounakis, S.~Greene, R.~Huang
  \emph{et~al.}, ``The power grid library for benchmarking ac optimal power
  flow algorithms,'' \emph{arXiv preprint arXiv:1908.02788}, 2019.

\bibitem{profile1}
S.~Kazarlis, A.~Bakirtzis, and V.~Petridis, ``A genetic algorithm solution to
  the unit commitment problem,'' \emph{IEEE Transactions on Power Systems},
  vol.~11, no.~1, pp. 83--92, 1996.

\bibitem{profile2}
C.~Grigg, P.~Wong, P.~Albrecht, R.~Allan, M.~Bhavaraju, R.~Billinton, Q.~Chen,
  C.~Fong, S.~Haddad, S.~Kuruganty, W.~Li, R.~Mukerji, D.~Patton, N.~Rau,
  D.~Reppen, A.~Schneider, M.~Shahidehpour, and C.~Singh, ``The ieee
  reliability test system-1996. a report prepared by the reliability test
  system task force of the application of probability methods subcommittee,''
  \emph{IEEE Transactions on Power Systems}, vol.~14, no.~3, pp. 1010--1020,
  1999.

\end{thebibliography}
\end{document}